\documentclass[12pt,final]{article}
\usepackage{amsmath,amsthm,amsfonts,amssymb,graphicx,enumerate,psfrag}
\usepackage{tikz,pgfplots}
\usepackage{mathrsfs}
\usepackage{fullpage}

\usepackage[colorlinks,citecolor=blue,urlcolor=blue]{hyperref}
\definecolor{myblue}{RGB}{51,51,178}
\definecolor{myred}{RGB}{189,26,26}
\definecolor{mygreen}{RGB}{0,128,0}
\usepackage[utf8]{inputenc} 
\newtheorem{lemma}{Lemma}

\newtheorem{theorem}{Theorem}
\newtheorem{example}{Example}
\newtheorem{proposition}{Proposition}

\newtheorem{remark}{Remark}

\newcommand{\dN}{\mathbb {N}}
\newcommand{\dZ}{\mathbb {Z}}

\newcommand{\dR}{\mathbb {R}}
\newcommand{\dd}{\mathrm{d}}

\newcommand{\EE}{{\mathbb{E}}}
\newcommand{\PP}{{\mathbb{P}}}
\newcommand{\dX}{{\mathbb{X}}}
\newcommand{\dY}{{\mathbb{Y}}}
\newcommand{\Lip}{{\mathrm{Lip}}}
\newcommand{\diam}{{\mathrm{diam}}}
\newcommand{\dist}{{\mathrm{dist}}}
\newcommand{\tmix}{\mathrm{t}_{\textsc{mix}}}

\newcommand{\tent}{\mathrm{t}_{\textsc{ent}}}
\newcommand{\tmls}{{\mathrm{t}_{\textsc{mls}}}}

\newcommand{\tv}{{\textsc{tv}}}
\newcommand{\Ent}{{\mathrm{Ent}}}

\newcommand{\Var}{{\mathrm{Var}}}

\newcommand{\Varent}{{\mathrm{Varent}}}

\title{A new cutoff criterion for  non-negatively curved chains}
\author{Francesco Pedrotti and Justin Salez}
\begin{document}
\maketitle
\begin{abstract}The cutoff phenomenon was recently shown to systematically follow from non-negative curvature and the product condition, for all Markov diffusions. The proof crucially relied on a classical \emph{chain rule} satisfied by the carré du champ operator, which is specific to differential generators and hence fails on discrete spaces. In the present paper, we show that an approximate version of this chain rule in fact always holds, with an extra cost that depends on the log-Lipschitz regularity of the considered observable. As a consequence, we derive a new cutoff criterion for non-negatively curved chains on finite spaces. The latter allows us to recover, in a simple and unified way, a number of historical instances of cutoff that had been established through   model-specific arguments. Emblematic examples include random walk on the hypercube, random transpositions, random walk on the multislice, or MCMC samplers for popular spin systems such as the Ising and Hard-core models on bounded-degree graphs.
\end{abstract}
\tableofcontents
\section{Introduction}
 \subsection{The cutoff phenomenon} 
\begin{figure}
\begin{center}
\pgfplotsset{compat=1.16}
\pgfplotsset{ticks=none}
\begin{tikzpicture}
	\begin{axis}[
		width = 14.5cm,
		height = 9cm,
		axis x line=middle,
		axis y line=middle,
		xlabel = $t$,
	    xlabel style ={at={(1,-0.1)}},
		ylabel style ={at={(0,1.05)}},
		ylabel = $\tv(X_t)$,
		clip = false,
		grid=both,
		grid style={dashed, line width=.5pt, draw=gray!10},
		xmode = normal,
		ymode = normal,
		line width = 1pt,
		legend cell align = left,
		legend style = {fill=none, at={(0.9,0.9)}, anchor = north east},
		yticklabel style={above left},
		anchor = north west,
		tickwidth={5pt},
		xtick align = outside,
		ytick align = outside,
		ymax = 1.05,
		ymin = 0,
		xmin = 0,
		xmax = 40,
		x axis line style=-,
		y axis line style=-,
		domain = 0:40,
		samples = 1000
		]
		\def\scale{3}
\addplot[solid,  myred, line width = 1.5 pt] {
-rad(atan(x-20))/rad(atan(20))/2 + 0.5	
	};
%
\draw[dashed, color= myblue] (15,0) -- (15,0.9515);
\draw[dashed, color= myblue] (0,0.9515) -- (15,0.9515);
\draw[dashed, color= mygreen] (25,0) -- (25,0.0484723);
\draw[dashed, color= mygreen] (0,0.0484723) -- (25,0.0484723);	

\draw[color= myblue] (0,0.9515) -- (-0.5,0.9515) node[color = black, anchor=east] {$\textcolor{myblue}{1-\varepsilon}$};
\draw[color= myblue] (15,0) -- (15,-0.02) node[color=black, anchor=north] {$\textcolor{myblue}{\tmix(1-\varepsilon)}$};
\draw[color= mygreen] (25,0) -- (25,-0.02) node[color = black, anchor=north] {${\textcolor{mygreen}{\tmix(\varepsilon)}}$};
\draw[color= mygreen] (0,0.0484723) -- (-0.5,0.0484723) node[color = black,anchor=east] {$\textcolor{mygreen}{\varepsilon}$};

\draw[color= black] (0,1.0) -- (-0.5,1.0) node[color = black,anchor=south east] {${1}$};

\end{axis}
	
\end{tikzpicture}
\caption{A typical plot of the distance to equilibrium $t\mapsto \tv(X_t)$. As the ratio $\frac{\tmix(1-\varepsilon)}{\tmix(\varepsilon)}$ approaches $1$, the transition to equilibrium becomes abrupt (cutoff).}
\label{fig:cutoff}
\end{center}
\end{figure}
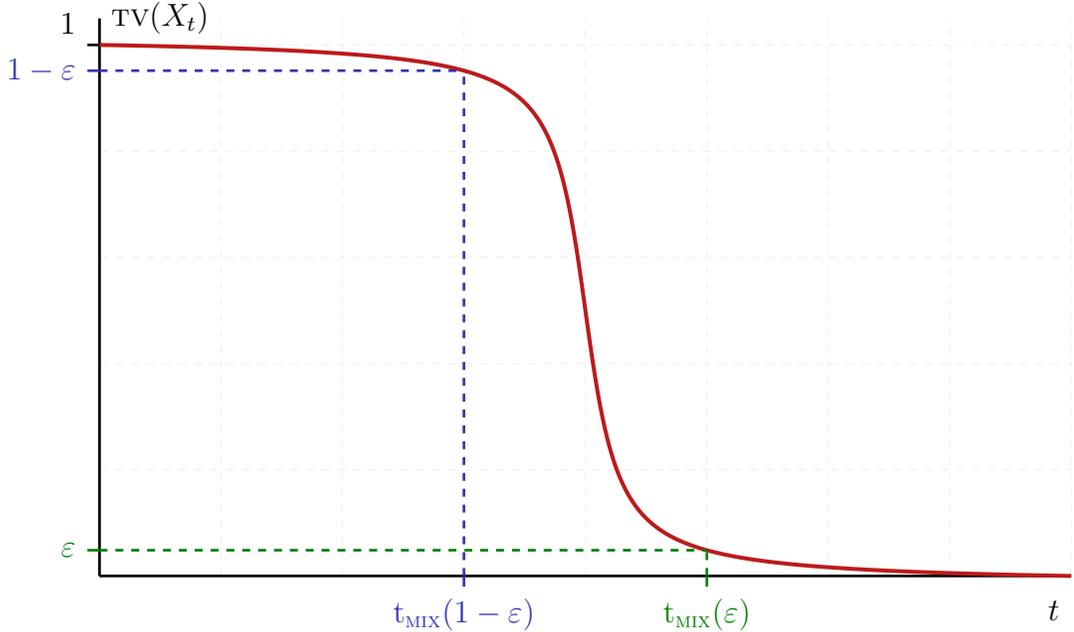

Consider a continuous-time Markov chain $(X_t)_{t\ge 0}$  on a finite state space  $\dX$.  Under the usual irreducibility assumption, the law of $X_t$ approaches a unique stationary distribution $\pi$ as $t\to\infty$, and it is natural to ask for the time-scale on which this convergence occurs. This is  formalized by the notion of \emph{mixing times} \cite{MR3726904}, defined for any precision $\varepsilon\in(0,1)$ by
\begin{eqnarray}
\label{def:tmix}
\tmix(\varepsilon) \ := \ \min\{t\ge 0\colon \tv(X_t)\le \varepsilon\}, & \textrm{where} &  \tv(X) \ = \ \sup_{A\subseteq\dX}\left|\PP(X\in A)-\pi(A)\right|.
\end{eqnarray}
In practice, the model under consideration often involves a natural size parameter $n\in\dN$ -- which will here remain implicit for notational ease -- and the interest is in the large-size limit $n\to\infty$. In certain cases, an abrupt transition  from out-of-equilibrium to equilibrium  has been observed, whereby  the distance to equilibrium $t\mapsto\tv(X_t)$ approaches a step function as $n\to\infty$, as illustrated on Figure \ref{fig:cutoff}. In other words, for any fixed precision $\varepsilon\in(0,1)$,
\begin{eqnarray*}
\frac{\tmix(1-\varepsilon)}{\tmix(\varepsilon)} & \xrightarrow[n\to\infty]{} & 1.
\end{eqnarray*}
This is the celebrated \emph{cutoff phenomenon}, discovered four decades ago in the context of card shuffling \cite{aldous1986shuffling,aldous1983mixing,diaconis1996cutoff}, and established since then in nearly a hundred different Markov chains arising in a broad variety of settings. Despite the accumulation of many examples, this phenomenon is
still far from being understood, and identifying the general conditions that trigger it has become one of the biggest challenges in the quantitative analysis of ergodic Markov processes. We refer the interested reader to the recent paper \cite{MR4780485} and the references therein for a detailed account of this fascinating question. In the present work, we provide a new, simple and unifying cutoff criterion for Markov chains that have non-negative curvature, in a sense that we now recall.

\subsection{Bakry-\'Emery curvature}
Introduced four decades ago in the context of diffusions on manifolds \cite{MR889476}, the Bakry-\'Emery theory of curvature is a powerful framework for the quantitative study of Markov semi-groups. We refer the unfamiliar reader to the  textbook \cite{MR3155209} for a comprehensive introduction. To keep the exposition simple, we shall here restrict our attention to finite state spaces, as considered, e.g., in \cite{MR1665591,MR3492631}. 
Upon rescaling time by a constant factor if needed, we may then assume  that the generator of our Markov process acts on functions $f\colon\dX\to\dR$ as follows:
\begin{eqnarray}
\label{def:L}
(Lf)(x) & = & \sum_{y\in\dX}T(x,y)\left(f(y)-f(x)\right),
\end{eqnarray}
for some stochastic matrix $T$ on $\dX$. The \emph{carré du champ} operator is obtained by squaring the discrete gradient in the above definition, and dividing by $2$:
\begin{eqnarray*}
(\Gamma f)(x) & = & \frac{1}{2}\sum_{y\in\dX}T(x,y)\left(f(y)-f(x)\right)^2.
 \end{eqnarray*} 
 Following Bakry and \'Emery \cite{MR889476}, we say that the chain is \emph{non-negatively curved} when the carré du champ operator sub-commutes with the semi-group $P_t=e^{tL}$, i.e.
\begin{eqnarray}
\label{assume:curved}
\forall t\ge 0,\qquad \Gamma P_t & \le & P_t\Gamma.
\end{eqnarray}
Thanks to the semi-group property $P_{t+s}=P_tP_s$, it is in fact enough to verify this when $t$ is infinitesimally small, leading to the more effective criterion $\Gamma_2\ge 0$ where $\Gamma_2$ denotes the iterated carré du champ operator. As a consequence, non-negative curvature is easy to check in practice, and several fundamental examples will be given in Section \ref{sec:app}. 
\subsection{Main result}
In addition to the  curvature condition (\ref{assume:curved}), we will require that the natural adjacency relation defined on our state space by $x\sim y\Longleftrightarrow T(x,y)>0$ is symmetric, i.e.
\begin{eqnarray}
\label{assume:symmetry}
\forall x,y\in\dX, \qquad T(x,y)>0 & \Longrightarrow & T(y,x)>0.
\end{eqnarray}
Note that this is weaker than the usual reversibility requirement $\pi(x)T(x,y)=\pi(y)T(y,x)$, which expresses the self-adjointness property $L^\star=L$ in $L^2(\pi)$. Let us define the \emph{degree} of the chain to be the inverse of the minimum non-zero transition probability:
\begin{eqnarray*}
d & := & \max_{x\sim y}\left\{\frac{1}{T(x,y)}\right\}.
\end{eqnarray*} 
This parameter controls the \emph{sparsity} of the transition matrix $T$: indeed, no row or column can have more than $d$ non-zero entries, and  $d$ is exactly the maximum degree of the graph in the special case of simple random walks. Next, we recall that the \emph{inverse modified log-Sobolev constant} is   the smallest number $\tmls$ such that the entropy decay
\begin{eqnarray}
\label{def:MLSI}
\forall t\ge 0,\qquad \Ent(X_t) & \le &  \Ent(X_0)\exp\left\{-\frac{t}{\tmls}\right\},
\end{eqnarray} 
holds for any initial condition $X_0$, where $\Ent(X)$ denotes the relative entropy of $X$ with respect to equilibrium, as defined at (\ref{def:ent}) below. Thanks to the semi-group property, it is  here again enough to consider the regime where $t$ is infinitesimally small, leading to a more effective variational characterization of $\tmls$ in terms of the Dirichlet form \cite{MR2283379,MR2341319}. Finally, let us be more explicit about the type of initial conditions that we allow. In the traditional literature on mixing times (see, e.g., \cite{MR3726904}), the Markov chain under consideration either starts from a designated ``origin'' $o\in\dX$, or from the ``worst'' possible initial distribution. To encompass both settings, we  will here consider mixing times of the form 
\begin{eqnarray*}
\tmix^{(S)}(\varepsilon) & := & \max_{o\in S}\left\{\tmix^{(o)}(\varepsilon)\right\},
\end{eqnarray*}
where  $S\subseteq\dX$ is an arbitrary (non-empty) region of allowed initial positions, and where the notation $\tmix^{(o)}(\varepsilon)$ refers to the particular initialization $X_0=o$. Note that by convexity of total variation,   $\tmix^{(S)}(\varepsilon)$ is in fact the worst-case mixing time over all initial distributions that are supported on $S$. The two standard settings mentioned above correspond to the extremal choices $S=\{o\}$ and $S=\dX$, respectively. However, we emphasize that our criterion below applies to \emph{any} region $S\subseteq\dX$. We are now ready to state our main result, in which the input data $(\dX,T,S)$ is simply referred to as a \emph{Markov triple}, and is assumed to depend on a parameter $n\in\dN$ which we keep implicit in the notation $\tmix,\tmls,d$. 

\begin{theorem}[Main result]\label{th:main}
Consider a sequence of Markov triples satisfying (\ref{assume:curved})-(\ref{assume:symmetry}) and
\begin{eqnarray*}
\frac{\tmix^{(S)}(\varepsilon)}{\tmls\log\log d} & \xrightarrow[n\to\infty]{} & +\infty,
\end{eqnarray*}
for some $\varepsilon\in(0,1)$. Then a cutoff occurs, i.e. for all $\varepsilon\in(0,1)$,
\begin{eqnarray*}
\frac{\tmix^{(S)}(1-\varepsilon)}{\tmix^{(S)}(\varepsilon)} & \xrightarrow[n\to\infty]{} & 1.
\end{eqnarray*}
\end{theorem}

\begin{remark}[Entropy mixing]As the careful reader will notice, our mixing-time upper bound is based on the modified log-Sobolev constant, and therefore controls mixing in the stronger entropy sense. 
As a consequence,  under the same assumptions, our proof actually yields 
\begin{eqnarray*}
\frac{\tmix^{(S)}(\varepsilon)}{\tent^{(S)}(\delta)} & \xrightarrow[n\to\infty]{} & 1,
\end{eqnarray*}
for any fixed $\varepsilon\in(0,1)$ and $\delta\in(0,\infty)$, where $\tent$ is the entropic mixing time obtained by replacing $\tv(X_t)$ by $\Ent(X_t)$ in the definition (\ref{def:tmix}). In words, cutoff occurs both in total-variation and relative entropy, at the same time.
\end{remark}

\section{Applications}
\label{sec:app}
Before diving into the proof of Theorem \ref{th:main}, let us demonstrate the effectiveness of our criterion by verifying it in a variety of historical examples where cutoff had been established through a delicate and model-specific analysis. We emphasize that the novelty here does not lie in the  results themselves, but rather in the unified and effortless way in which we recover them. 
\subsection{Conjugacy-invariant random walks on groups}
Consider a finite group $\dX$, equipped with a probability measure $\mu$ whose support is symmetric and generates the group. By definition, the (left) random walk on $\dX$ with increment law  $\mu$ is  the Markov chain on $\dX$ with transition matrix 
\begin{eqnarray*}
T(x,y) & := & \mu(yx^{-1}).
\end{eqnarray*}
By symmetry, the choice of the initial state is irrelevant, and we take it to be the identity element. In this context, the curvature assumption (\ref{assume:curved}) is well known to hold as soon as
\begin{eqnarray*}
\forall x,y\in\dX,\qquad \mu(xy) & = & \mu(yx). 
\end{eqnarray*}
We refer the interested reader to \cite{MR1665591,hermon2024concentrationinformationdiscretegroups} for a proof. Note that this property is trivially satisfied, in particular, when the group $\dX$ is Abelian. The simplest example is of course simple random walk on the boolean hypercube, which is well known to exhibit  cutoff. 
\begin{example}[Random walk on the hypercube]\label{ex:cube}Let $\dX$ be the additive group $\{0,1\}^n$, and $\mu$ the uniform distribution on its canonical basis. Then, 
\begin{equation*}
d \ = \ n, \qquad \tmls=\Theta(n),\qquad \tmix \ =\ \Theta(n\log n),
\end{equation*}
see \cite[Example 3.7]{MR2283379}. Thus, our criterion is satisfied and cutoff follows. 
\end{example}
We next consider the non-Abelian case of random transpositions, for which the occurrence of a cutoff  is a celebrated historical result due to Diaconis and Shahshahani \cite{diaconis1981generating}.
\begin{example}[Random transpositions]\label{ex:transpositions}Let $\dX$ be the symmetric group of order $n$, and $\mu$ the uniform measure on the set of transpositions. Then, the corresponding random walk satisfies 
\begin{equation*}
d \ = \ \Theta(n^2), \qquad \tmls=\Theta(n),\qquad \tmix \ =\ \Theta(n\log n),
\end{equation*}
see \cite{MR2023890,MR2094147}. Thus, our  criterion is again satisfied and cutoff follows.  
\end{example}
More generally, one can replace the set of transpositions in the above example by any conjugacy class whose \emph{complexity} (number of non-fixed points) is not too large.
\begin{example}[Random walks generated by a conjugacy class]Let $\dX$ be the symmetric group of order $n$, and $\mu$ the uniform measure on a non-trivial, symmetric conjugacy class $S\subseteq|\dX|$. Let  $k$ denote the number of non-fixed points in any member of  $S$. Then, 
\begin{equation*}
d \ = \ \Theta(n^k), \qquad \tmls=\Theta\left(\frac{n}{k}\right),\qquad \tmix \ =\ \Theta\left(\frac{n\log n}{k}\right),
\end{equation*}
so that our criterion is satisfied as long as $k=n^{o(1)}$, see again \cite{MR2094147}.  Note that the previous example corresponds to the special case where $k=2$. Interestingly,  cutoff is known to occur in the more general regime where $k=o(n)$, as conjectured by   Diaconis and Shahshahani \cite{diaconis1981generating}, and recently proved by Berestycki and \c{S}eng\"{u}l  \cite{MR3936154}. 
\end{example}
\subsection{Markovian projections}
An elementary (but seemingly new) observation about the curvature condition (\ref{assume:curved}) is that it is preserved under projections, in the following sense. Consider a surjective map $\Phi\colon\dX\to\dY$ from our state space onto another one. It is well known that the image of our Markov chain $(X_t)_{t\ge 0}$ under $\Phi$ is again a Markov chain, provided that for each $y\in\dY$, the quantity
$
\sum_{z\in\Phi^{-1}(y)}T(x,z)
$
depends on the state $x$ only through  $\Phi(x)$. In other words, we can write 
\begin{eqnarray*}
\sum_{z\in\Phi^{-1}(y)}T(x,z) & = & \widehat{T}\left(\Phi(x),y\right),
\end{eqnarray*}
for some matrix $\widehat{T}\colon\dY^2\to[0,1]$. The latter is then necessarily stochastic, and  it is nothing but the transition matrix of the Markovian projection  $\left(\Phi(X_t)\right)_{t\ge 0}$. By linearity, the above relation implies the identity $T(f\circ\Phi)=(\widehat{T}f)\circ\Phi$ for all observables $f\colon\dY\to\dR$, and a similar intertwining relation holds at the level of semi-groups and carré du champ operators: 
\begin{eqnarray*}
P_t(f\circ\Phi) \ = \  (\widehat P_t f)\circ\Phi,\quad \textrm{and} \quad
\Gamma(f\circ\Phi) \ = \  (\widehat \Gamma f)\circ\Phi.
\end{eqnarray*}
In particular, it readily follows from those identities that the sub-commutation relation $\Gamma P_t\le P_t\Gamma$ implies $\widehat{\Gamma}\widehat{P_t}\le \widehat{P_t}\widehat{\Gamma}$. Let us record this fact for future reference.
\begin{lemma}[Non-negative curvature is preserved under projection]Any Markovian projection of a non-negatively curved Markov chain is again non-negatively curved. 
\end{lemma}
Similarly, many important theoretical parameters of Markov chains can only improve under Markovian projections. This classically includes   the inverse modified log-Sobolev constant $\tmls$ and the degree $d$, making our theorem particularly well-behaved under projections. Let us illustrate this general principle with a few emblematic examples. 
\begin{example}[Ehrenfest model]Consider $n$ unlabeled particles evolving between two containers as follows: at unit rate, a particle is chosen uniformly at random, and moved from its container to the other one. This simple model of diffusion was famously proposed by Tatiana and Paul Ehrenfest to explain the second law of thermodynamics. Formally, it can be obtained by projecting the random walk on the hypercube (Example \ref{ex:cube}) through the function
\begin{eqnarray*}
\Phi(x_1,\ldots,x_n) & = & x_1+\cdots+x_n.
\end{eqnarray*}
An easy application of Wilson's method \cite{MR2023023} shows that the mixing time of this process remains of order $n\log n$, as for random walk on the hypercube. Since the parameters $d$ and $\tmls$ can only decrease, the cutoff criterion which was verified on the hypercube remains satisfied here, and we recover a celebrated historical result, proved in the seminal lecture notes \cite{aldous1983mixing}.
\end{example}

\begin{example}[Bernoulli-Laplace model]Consider $n$ unlabeled balls, half of which are red and placed in a first urn, the other half being blue and placed in a second urn.  A transition consists in choosing a pair of balls  at random and swapping their positions. The  system is clearly determined by the number of red balls in the first urn, and its evolution can be realized by projecting the transposition walk (Example \ref{ex:transpositions}) through the map
\begin{eqnarray*}
\Phi(\sigma) & := & \sum_{i=1}^{n/2}{\bf 1}_{\sigma(i)\le n/2}.
\end{eqnarray*}
The mixing time is easily seen to be at least of order $n\log n$, just like in Example \ref{ex:transpositions}. Thus, our cutoff criterion is satisfied and cutoff follows, as first shown in \cite{diaconis1981generating}.
\end{example}

\begin{example}[Random walk on the multislice] Fix a finite sequence of positive integers $\kappa=(\kappa_1,\ldots,\kappa_L)$, and consider the set of words of length $n:=\kappa_1+\cdots+\kappa_L$ in which each symbol $\ell\in[L]$ appears exactly $\kappa_\ell$ times:
\begin{eqnarray*}
\Omega_\kappa & := & \left\{\omega=(\omega_1,\ldots,\omega_n)\in[L]^n\colon \sum_{i=1}^n{\bf 1}_{(\omega_i=\ell)} = \kappa_\ell \textrm{ for each }\ell\in[L] \right\}.
\end{eqnarray*}
This natural combinatorial set is sometimes called a \emph{multislice}, and the random walk on it consists in  swapping two uniformly chosen coordinates at rate $1$; see \cite{MR4353960} and the references therein. This is the image of  the transposition walk (Example \ref{ex:transpositions}) through the map
\begin{eqnarray*}
\Phi(\sigma) \ := \ \left(\phi(\sigma_1),\ldots,\phi(\sigma_n)\right), & \textrm{ where } & \phi(i) = \min\left\{\ell\le L\colon \kappa_1+\cdots+\kappa_\ell\ge i\right\}.
\end{eqnarray*}
The parameters $d$ and $\tmls$ can only decrease compared to their values in Example \ref{ex:transpositions}. On the other hand,  $\tmix$ is easily seen to be at least  of order $n\log\left(n-\kappa_{\max}\right)$, where $\kappa_{\max}:=\max_\ell\kappa_\ell$. Thus, our criterion predicts a cutoff as soon as $\frac{\log\left(n-\kappa_{\max}\right)}{\log\log n}\to+\infty$. In particular, we can take $\kappa_1=\cdots=\kappa_L=n/L$ where $L\ge 2$ is fixed to recover the main result of \cite{MR1436486}.
\end{example}
\subsection{MCMC samplers}
Consider a fully supported probability measure $\pi$ on our finite state space $\dX$. 
The celebrated \emph{Markov chain Monte Carlo revolution} in computational statistics is fundamentally based on the simple but far-reaching idea -- attributed to Metropolis \cite{doi:10.1063/1.1699114} and Hastings \cite{MR3363437} -- that approximate samples from $\pi$ can be efficiently produced by running an appropriate Markov chain that admits $\pi$ has its equilibrium law; see the survey paper by P. Diaconis \cite{MR2476411} and the references therein. Following \cite{MR3646066,MR4414692,pedrotti2023contractivecouplingratescurvature}, we will  here focus on implementations of the form
\begin{eqnarray}
\label{def:glauber}
(Lf)(x) & := & \frac{1}{|G|}\sum_{\tau\in G}\sqrt{\frac{\pi(\tau x)}{M\pi(x)}}\left(f(\tau x)-f(x)\right),
\end{eqnarray}
where $G$ is a given set of maps $\tau\colon x\mapsto \tau x$ on $\dX$ describing the allowed moves, and where
\begin{eqnarray*}
M & := & \max_{x\in\dX,\tau\in G}\left\{\frac{\pi(\tau x)}{\pi(x)}\right\},
\end{eqnarray*}
is a normalizing constant which is irrelevant for cutoff but ensures that the chain jumps at rate at most $1$, in compliance with our  convention  (\ref{def:L}). The dynamics (\ref{def:glauber}) is clearly reversible w.r.t. the target measure $\pi$. To gain some intuition, consider the instructive case where $\pi$ is the uniform measure on $\dX=\{0,1\}^n$ and $G=\{\tau_1,\ldots,\tau_n\}$, where $\tau_i\colon\dX\to\dX$ is  the map that flips the $i-$th coordinate. In that case, (\ref{def:glauber})  is exactly the generator of the random walk on the hypercube (Example \ref{ex:cube}), which has been seen to exhibit cutoff. In light of this, it is natural to expect a similar phenomenon when sampling from more general high-dimensional measures with weak dependencies. In an impressive series of works \cite{MR3020173,lubetzky2014cutoff,MR3434254,MR3486171,MR3729612}, Lubetzky and  Sly developed a very sophisticated framework named \emph{Information Percolation}, which enabled them to confirm the above intuition for various high-temperature spin systems on arbitrary bounded-degree graphs, such as the celebrated Ising and Hard-core models. As we will now see, our main criterion is easily verified in those emblematic models. 

\begin{example}[Ising model on a graph]The Ising model with inverse temperature $\beta\ge 0$ on a finite graph $\mathbb G=(\mathbb V,\mathbb E)$ is the probability measure
\begin{eqnarray*}
\pi(x) \ \propto \ \exp\left\{\beta \sum_{\{i,j\}\in \mathbb  E}x_ix_j\right\} & \textrm{ on } &  \dX=\{-1,1\}^{\mathbb V}.
\end{eqnarray*}
Consider the sampler (\ref{def:glauber}) with allowed moves $G=\{\tau_i\}_{i\in \mathbb V}$, where $\tau_i\colon\dX\to\dX$ is  the map that flips the $i-$th coordinate. Writing $\Delta$ for the maximum degree in  $\mathbb G$, we have $M\le e^{2\beta\Delta}$ and $d\le |\mathbb V|e^{2\beta\Delta}$. Now, it follows from \cite{pedrotti2023contractivecouplingratescurvature} that this chain is non-negatively curved whenever
\begin{eqnarray}
\label{assume:beta}
\Delta (1-e^{-2\beta})e^{2\Delta \beta} & \le & 1.
\end{eqnarray}
(The result therein is stated in the context where $\mathbb G$ is a  subgraph of $\dZ^d$, but the proof never uses this).  Moreover, under this condition, the same work (or \cite{MR3646066,MR4414692}) gives $\tmls\le |\mathbb V|\sqrt{M}\le |\mathbb V|e^{\beta\Delta}$.  Since $\tmix$ is at least of order $|\mathbb V|\log |\mathbb V|$, our criterion holds along any sequence of bounded-degree graphs with diverging size, in the high-temperature regime (\ref{assume:beta}).
\end{example}

\begin{example}[Hard-core model on a graph]The Hard-core model with fugacity $\lambda\in(0,1)$ on a finite graph $\mathbb G=(\mathbb V,\mathbb E)$ is the probability measure 
\begin{eqnarray*}
\pi(x)  \propto \lambda^{\sum_{i\in\mathbb V}x_i}
& \textrm{ on } & 
\dX \ = \left\{x\in \{0,1\}^{\mathbb V}\colon \forall \{i,j\}\in\EE, x_ix_j=0\right\}.
\end{eqnarray*}
Consider the associated sampler (\ref{def:glauber}) with  $G=\{\tau_i\}_{i\in \mathbb V}$, where $\tau_i\colon\dX\to\dX$ is  the map that flips the $i-$th coordinate if the resulting vector is in $\dX$, and does nothing otherwise. Note that $M=\lambda^{-1}$ and $d=|\mathbb V|{\lambda^{-1}}$. Again, non-negative curvature follows from \cite{pedrotti2023contractivecouplingratescurvature} as soon as 
\begin{eqnarray}
\label{assume:lambda}
\lambda\Delta & \le & 1,
\end{eqnarray}
where $\Delta$ denotes the maximum degree in  $\mathbb G$.  Moreover, under this condition, the very same work, or \cite{MR4414692}, implies that $\tmls\le |\mathbb V|\sqrt{M\lambda^{-1}}=|\mathbb V|\lambda^{-1}$. Since $\tmix$ is at least of order $|\mathbb V|\log |\mathbb V|$, our criterion holds along any sequence of bounded-degree graphs with diverging sizes, throughout the low-fugacity regime (\ref{assume:lambda}).
\end{example}

\section{Proof}
Following the ideas exposed  in \cite{MR4780485,MR4765357,hermon2024concentrationinformationdiscretegroups,salez2025cutoffnonnegativelycurveddiffusions}, we will estimate the width of the mixing window through the information-theoretic notions of \emph{entropy} and \emph{varentropy}. Recall that those statistics are respectively defined, for any $\dX-$valued random variable $X$, as 
\begin{eqnarray}
\label{def:ent}
\Ent(X) \ := \ \EE[\log f(X)] & \textrm{ and } & \Varent(X) \ := \ \Var\left[\log f(X)\right],
\end{eqnarray} 
where $f$ is the density of $X$ w.r.t. the stationary measure $\pi$.
The key ingredient in the recent breakthrough \cite{salez2025cutoffnonnegativelycurveddiffusions}  was the observation that the entropy and varentropy of a non-negatively curved diffusion $(X_t)_{t\ge 0}$ are related through the  differential inequality
\begin{eqnarray}
\label{assume:diff}
\frac{\dd}{\dd t}\,\Ent(X_t) & \le & -\frac{\Varent(X_t)}{2t},
\end{eqnarray}
for all $t\ge 0$. Unfortunately, this crucially relied on the so-called \emph{chain rule} satisfied by the associated carré du champ operator, which notoriously fails in the discrete setup considered here. Nevertheless, we will now show that an approximate version of the chain rule actually always holds, with a multiplicative error that depends on the log-Lipschitz regularity of the considered observable. This is reminiscent of the regularization principle used in \cite{MR4765699,MR4620718} to relate the log-Sobolev constant and its modified version.
\subsection{Approximate chain rule}\label{sec:chainrule} In the very different context of Markov diffusions  on the $d-$dimensional Euclidean space, the generator and carré du champ operator act on smooth functions $f\colon\dR^d\to\dR$ as follows:
\begin{eqnarray*}
Lf \ = \ \sum_{i=1}^dg_i \partial_i f+ \sum_{i,j=1}^d g_{ij}\partial_{ij}f & \textrm{ and } & 
\Gamma f\ = \ \sum_{i,j=1}^d g_{ij}\partial_i f \partial_j f,
\end{eqnarray*}
for a given collection of functions $(g_i)_{1\le i\le d}$ and $(g_{ij})_{1\le i,j\le d}$. It then easily follows that 
\begin{eqnarray}
\label{chain}
\Gamma(\log f) & = & \frac{Lf}{f}-L\log f,
\end{eqnarray}
for any smooth positive function $f$. As explained above, this crucial chain rule notoriously fails on our discrete state space $\dX$. However, an approximate version of it turns out to hold, with an extra ``cost'' that accounts for the intrinsic roughness of $\log f$. Specifically, let
\begin{eqnarray*}
\Lip(f) & := & \max_{x\sim y}\left|f(x)-f(y)\right|,
\end{eqnarray*}
denote the Lipschitz constant of a function $f\colon\dX\to\dR$, and let us introduce the cost function
\begin{eqnarray*}
\Psi(r) & := & \frac{r^2}{2\left(r +e^{-r}-1\right)},
\end{eqnarray*}
with the understanding that $\Psi(0)=1$. This function is easily seen to be continuously increasing from $\Psi(-\infty)=0$ to $\Psi(+\infty)=+\infty$, and to satisfy $\Psi(r)\le 1+r$ for all $r\ge 0$.  We then have the following approximate chain rule,  of which the identity (\ref{chain}) can be seen as the infinitely-smooth limit $r=0$. 
\begin{lemma}[Approximate chain rule]\label{lm:chain}Fix $f\colon\dX\to(0,\infty)$ and set $r:=\Lip(\log f)$.  Then, 
\begin{eqnarray*}
\Psi(-r)\,\left(\frac{Lf}{f}-L\log f\right)\ \le & \Gamma(\log f) & \le \ \Psi(r)\,\left(\frac{Lf}{f}-L\log f\right).
\end{eqnarray*}
\end{lemma}
\begin{proof}Since $\Psi$ increases on $\dR$, we have $\Psi(-r)\le \Psi(\ell)\le\Psi(r)$ whenever $\ell\in[-r,r]$. In particular, if $x,y\in\dX$ are neighbors, we may take $\ell=\log\frac{f(x)}{f(y)}$ and $r:=\Lip(\log f)$ to get
\begin{equation*}
\Psi(-r)\left(\frac{f(y)}{f(x)}-1+\log\frac{f(x)}{f(y)}\right) \ \le \ \frac{1}{2}\log^2\frac{f(x)}{f(y)} \ \le \ \Psi(r)\left(\frac{f(y)}{f(x)}-1+\log\frac{f(x)}{f(y)}\right).
\end{equation*}
Multiplying through by $T(x,y)$ and summing over all $y\in\dX$ concludes the proof. 
\end{proof}
As promised, this approximate chain rule allows us to establish a version of the information differential inequality (\ref{assume:diff}) for all non-negatively curved Markov chains on finite spaces.
\begin{proposition}[Information-differential inequality]\label{pr:diff}Consider  a non-negatively curved Markov chain  $(X_t)_{t\ge 0}$ on a finite state space, starting from a deterministic point. Then,
\begin{eqnarray*}
\frac{\dd}{\dd t}\,\Ent(X_t) & \le & -\frac{\Varent(X_t)}{2t(1+\Lip(\log f_t))},
\end{eqnarray*}
where $f_t$ is the density of $X_t$ with respect to $\pi$.
\end{proposition}
\begin{proof}The sub-commutation relation (\ref{assume:curved}) classically provides the following \emph{local Poincaré inequality} along the chain (see, e.g., \cite{MR3155209}): for any observable $g\colon\dX\to\dR$ and any time $t\ge 0$, 
\begin{eqnarray*}
\Var\left[g(X_t)\right]  & \le & 2t\,\EE\left[\Gamma g(X_t)\right].
\end{eqnarray*}
Applying this to $g=\log f_t$ and using our chain rule together with $\Psi(r)\le 1+r$, we obtain
\begin{eqnarray*}
\Varent(X_t)  & \le &  2t\,\left(1+\Lip(\log f_t)\right)\,\EE\left[\left(\frac{Lf_t}{f_t}-L\log f_t\right)(X_t)\right].
\end{eqnarray*}
Now, since $X_t$ has law $f_t\,\dd\pi$, the expectation appearing on the right-hand side reads
\begin{eqnarray*}
\EE\left[\left(\frac{Lf_t}{f_t}-L\log f_t\right)(X_t)\right] & = & \pi\left[Lf_t-f_tL\log f_t\right] \ =\ -\pi\left[f_tL\log f_t\right],
\end{eqnarray*}
because $\pi L=0$.
On the other hand, the Fokker-Planck equation $\frac{\dd f_t}{\dd t}=L^\star f_t$ shows that
\begin{eqnarray*}
\frac{\dd \Ent(X_t)}{\dd t} & = & \frac{\dd}{\dd t}\pi\left[f_t\log f_t\right] \ =\ \pi\left[(L^\star f_t)(1+\log f_t)\right] \ = \  \pi\left[f_tL\log f_t\right],
\end{eqnarray*}
where we have used the very definition of the adjoint operator $L^\star$ in $L^2(\pi)$, and the mass conservation property $L1=0$.  Combining the last three displays concludes the proof.
\end{proof}
\subsection{Spatial regularity of the information content}
To turn Proposition \ref{pr:diff} into an effective statement, we  need to  estimate  the new regularity term $\Lip(\log f_t)$ featuring in our information-differential inequality. This is the content of the following lemma, in which $\diam$ denotes the diameter of the chain, or more accurately, of the graph induced by the adjacency relation $\sim$. 
\begin{lemma}[Spatial regularity of the heat kernel]\label{lm:reg}Let $(X_t)_{t\ge 0}$ be a Markov chain satisfying the symmetry (\ref{assume:symmetry}), and let $f_t$ denote the density of $X_t$ w.r.t. equilibrium. Then, for all $t\ge 0$, 
\begin{eqnarray*}
\Lip\left(\log f_t\right) & \le & 3+3\log d+3\log\left(1\vee\frac{\diam}{4t}\right),
\end{eqnarray*}
where $d$ and $\diam$ denote the degree and diameter of the chain, respectively.
\end{lemma}
\begin{proof}
Clearly, $f_t$ is a convex combination of the extremal densities $f_t^{(o)},o\in\dX$, where $f_t^{(o)}$ denotes the density of $X_t$ in the special case where $X_0=o$. As a consequence, it is enough to prove the result when $X_0$ is deterministic. This was actually done in \cite[Lemma 10]{MR4780485}, but only in the regime where $t\ge \diam/4$, in which case the last term simply vanishes. Since  this was valid for any transition matrix $T$ with symmetric support, we may fix $\theta\in(0,1]$ and apply it to the modified transition matrix
\begin{eqnarray*}
\widehat{T} & := & \theta T+(1-\theta)\mathrm{Id}.
\end{eqnarray*}
Note that this transformation preserves the adjacency relation $\sim$, hence the diameter. On the other hand, the associated semi-group becomes $ \widehat{P}_t  =  P_{\theta t}$, while  the minimum non-zero transition probability satisfies $
  \widehat{d} \le  d/\theta$. Thus, the conclusion now reads
  \begin{eqnarray*}
 \Lip\left(\log f_t\right) & \le & 3+3\log\frac{d}{\theta},
 \end{eqnarray*} 
provided that the constraint $t/\theta\ge \diam/4$ is satisfied. Since this is true for any choice of $\theta\in(0,1]$, we may finally optimize the bound by choosing $\theta:=\min\left\{1,\frac{4t}{\diam}\right\}$.
\end{proof}
Recalling that our target criterion for cutoff only involves the parameters $d$ and $\tmls$, we would now like to estimate the diameter appearing in the above lemma in terms of $d$ and $\tmls$. This is the content of the following lemma, which appears to be new.
\begin{lemma}[Diameter and modified log-Sobolev constant]\label{lm:diam}We always have
\begin{eqnarray*}
\diam & \le & 16\,\tmls\log{2d}.
\end{eqnarray*}
\end{lemma}
\begin{proof}
Let us first assume that our transition matrix $T$ is lazy and reversible. Then, by virtue of a classical argument due to Herbst (see, e.g., \cite{MR1849347}), the modified log-Sobolev inequality guarantees sub-Gaussian concentration under the stationary measure. More precisely, for any function $f\colon\dX\to\dR$ with $\pi[f]=0$ and $\Lip(f)\le 1$, and any $t\ge 0$, we have
\begin{eqnarray*}
\log \pi\left[e^{tf}\right] & \le & \frac{\tmls t^2}{4}. 
\end{eqnarray*}
Using the crude bound $\pi[e^{tf}]\ge \pi_{\min} e^{t\max f}$, we obtain 
\begin{eqnarray*}
t\max f & \le & \frac{\tmls t^2}{4}+\log\frac{1}{\pi_{\min}}. 
\end{eqnarray*}
This is valid for any $t\ge 0$, and the optimal choice $t=2\max f/\tmls$ yields
\begin{eqnarray*}
\max f & \le & \sqrt{\tmls\log\frac{1}{\pi_{\min}}}.
\end{eqnarray*} 
Of course, the same bound applies to $-f$, so we conclude that
\begin{eqnarray*}
\max f-\min f & \le & 2\sqrt{\tmls\log\frac{1}{\pi_{\min}}}.
\end{eqnarray*} 
Since this is invariant under shifting $f$ by a constant, our assumption $\pi[f]=0$ can now be dropped. In particular, we can take $f(x)=\dist(o,x)$, where $\dist(\cdot,\cdot)$ is the graph distance induced by the adjacency relation $\sim$. Since the base-point $o\in\dX$ is arbitrary, we obtain 
\begin{eqnarray*}
 \diam & \le & 2\sqrt{\tmls\log\frac{1}{\pi_{\min}}}.
 \end{eqnarray*} 
Finally, note that the matrix $Q:=T^\diam$ is stochastic and satisfies $\pi Q=\pi$ as well as $Q(x,y)\ge d^{-\diam}$ for all $x,y\in\dX$. Consequently, we have
\begin{eqnarray*}
\forall y\in\dX,\qquad \pi(y) & = & \sum_{x\in\dX}\pi(x)Q(x,y) \ \ge \ d^{-\diam},
\end{eqnarray*}
i.e. $\pi_{\min}\ge d^{-\diam}$. Inserting this into the previous display and simplifying yields
 \begin{eqnarray*}
 \diam & \le & 4\tmls\log d.
 \end{eqnarray*} 
Now, this was established under the extra assumption that $T$ is lazy and reversible. In the general case, we can always apply the above inequality to the lazy reversible matrix $\widehat{T}:=\frac{T+T^\star+2\mathrm{Id}}{4}$, which satisfies $\widehat{\diam}=\diam$, $\widehat{\tmls}\le 4\tmls$ and $\widehat{d}\le 2d$. 
\end{proof} 
\subsection{The information-differential route to cutoff}\label{sec:final}
With the estimates of Lemmas \ref{lm:reg}-\ref{lm:diam} at hand, our information-differential inequality (Proposition {\ref{pr:diff}}) becomes fully effective, and we may finally use it to deduce cutoff. To this end, let us recall that the relative entropy of a $\dX-$valued  variable $X$ always provides an upper-bound on its total variation distance to equilibrium, as per the celebrated Pinsker inequality: 
\begin{eqnarray}
\label{pinsker}
2\tv^2(X) & \le & \Ent(X).
\end{eqnarray}
Varentropy allows us to reverse this inequality, as established in  \cite[Lemma 8]{MR4780485}:
\begin{eqnarray}
\label{reverse}
\Ent(X) & \le & \frac{1+\sqrt{\Varent(X)}}{1-\tv(X)}.
\end{eqnarray}
This will play a crucial role in our proof. We will also use the standard mixing-time bound
\begin{eqnarray}
\label{tmixbound}
\tmix(\varepsilon) & \le & t+\tmls\,\log\left(1\vee \frac{\Ent(X_t)}{2\varepsilon^2}\right),
\end{eqnarray}
valid for any Markov chain $(X_t)_{t\ge 0}$ and any time $t\ge 0$, and which readily follows from (\ref{def:MLSI}), (\ref{pinsker}), and the semi-group property. We are now ready to prove Theorem \ref{th:main}.

\begin{proof}[Proof of Theorem \ref{th:main}]Consider a Markov triple $(\dX,T,S)$ as in Theorem \ref{th:main}. Fix $\varepsilon\in(0,1/2)$ and set $t_0:=\tmix^{(S)}(1-\varepsilon)$. Now, consider a continuous-time Markov chain $(X_t)_{t\ge 0}$ with transition matrix $T$ starting from a fixed state $o\in S$. Write $f_t$ for the density of $X_t$ w.r.t. equilibrium. First, Lemmas \ref{lm:reg}-\ref{lm:diam} and the inequality $\log u\le u-1$ ensure that for all $t>0$,
\begin{eqnarray*}
1+\Lip(\log f_t) & \le & 15\,\left(\log d+\frac{\tmls}{t}\right).
\end{eqnarray*}
On the other hand, for all $t\ge t_0$, we have $\tv(X_t)\le 1-\varepsilon$, so that  (\ref{reverse}) implies
\begin{eqnarray*}
\Varent(X_t) & \ge & (\varepsilon\Ent(X_t)-1)^2.
\end{eqnarray*}
In view of Proposition \ref{pr:diff}, we deduce that on $[t_0,\infty)$, we have the differential inequality
\begin{eqnarray*}
\frac{\dd}{\dd t}\,\Ent(X_t) & \le & -\frac{\left(\varepsilon\Ent(X_t)-1\right)^2}{30(t\log d+ \tmls)}.
\end{eqnarray*}
Integrating this inequality, and using $\log u\ge 1-\frac{1}{u}$, we obtain
\begin{eqnarray*}
\frac{1}{\varepsilon\Ent(X_t)-1} & \ge & \frac{1}{\varepsilon \Ent(X_{t_0})-1}+\frac{\varepsilon}{30\log d}\log\left(\frac{t\log d+\tmls}{t_0\log d+\tmls}\right)\\
& \ge &  \frac{\varepsilon(t-t_0)}{30(t\log d+\tmls)}.
\end{eqnarray*}
provided $t> t_0$ and $\Ent(X_t)>\frac{1}{\varepsilon}$. Consequently, for all $t> t_0$,
\begin{eqnarray*}
\Ent(X_t) & \le & \frac{1}{\varepsilon}+\frac{30(t\log d+\tmls)}{\varepsilon^2(t-t_0)}.
\end{eqnarray*}
Inserting this into the mixing-time bound (\ref{tmixbound}), we obtain
\begin{eqnarray*}
\tmix^{(o)}(\varepsilon) & \le & t+\tmls\log\left\{\frac{30}{\varepsilon^4}\right\}+\tmls\log\left\{1+\frac{t\log d+\tmls}{t-t_0}\right\}.
\end{eqnarray*}
This bound is valid for any $t>t_0$, so we may choose $t=t_0+\tmls$ to arrive at 
\begin{eqnarray*}
\tmix^{(o)}(\varepsilon) & \le & t_0+\tmls\log\left\{\frac{120}{\varepsilon^4}\right\}+\tmls\log\log d+\tmls\log\left\{2+\frac{t_0}{\tmls}\right\}.
\end{eqnarray*}
Taking a maximum over all $o\in S$ and recalling our choice for $t_0$, we conclude that
\begin{eqnarray*}
\label{window}
\tmix^{(S)}(\varepsilon) - \tmix^{(S)}(1-\varepsilon) & \le & \tmls\log\left\{\frac{120}{\varepsilon^4}\right\}+\tmls\log\log d+\tmls\log\left\{2+\frac{\tmix^{(S)}(1-\varepsilon)}{\tmls}\right\}.
\end{eqnarray*}
Finally, assume that our Markov triple $(\dX,T,S)$ depends on a parameter $n\ge 1$, and that
\begin{eqnarray*}
\frac{\tmix^{(S)}(\delta)}{\tmls\log\log d} & \xrightarrow[n\to\infty]{} & +\infty.
\end{eqnarray*}
for some fixed $\delta\in(0,1)$. Then, choosing $\varepsilon$ smaller than $1-\delta$ ensures that ${\tmix^{(S)}(1-\varepsilon)}\le\tmix^{(S)}(\delta)$, so that the last two displays together imply 
\begin{eqnarray*}
\frac{\tmix^{(S)}(\varepsilon) - \tmix^{(S)}(1-\varepsilon)}{\tmix^{(S)}(\delta)} & \xrightarrow[n\to\infty]{} & 0.
\end{eqnarray*}
Since this holds for arbitrarily small values of $\varepsilon$, cutoff follows.
\end{proof}
 \section*{Acknowledgment.}This work is supported by the ERC consolidator grant CUTOFF (101123174). Views and opinions expressed are however those of the authors only and do not necessarily reflect those of the European Union or the European Research Council Executive Agency. Neither the European Union nor the granting authority can be held responsible for them.

\bibliographystyle{plain}
\bibliography{draft}

\end{document}